\documentclass[12pt,a4paper]{article}
\usepackage{amsmath,amsthm,amssymb}
\usepackage[truedimen,margin=25truemm]{geometry}
\usepackage{mathrsfs}
\usepackage{here}
\usepackage{bm}
\usepackage{graphicx}
\usepackage[utf8]{inputenc}
\usepackage{float}

\theoremstyle{definition}
\newtheorem{definition}{定義}
\newtheorem{thm}[definition]{Theorem}

\newtheorem{cor}[definition]{Corollary}

\newtheorem{remark}[definition]{Remark}

\newcommand{\Mat}{{\rm Mat}}

\DeclareMathOperator{\SL}{SL}

\begin{document}

\title{Group Determinants and Invariant Rings}
\author{Yuka Yamaguchi and Naoya Yamaguchi}
\date{\today}

\maketitle

\begin{abstract}
In the study of group determinants, Frobenius introduced certain partial differential operators. 
This paper presents several results concerning the invariant rings derived from these partial differential operators. 
\end{abstract}

\section{Introduction}
Let $G = \{ g_{1}, g_{2}, \ldots, g_{n} \}$ be a finite group of order $n$, 
where $g_{1}$ is the identity element, 
and let $x_{g}$ be a variable associated with each $g \in G$. 
For an $m$-tuple of variables $\bm{z} := {}^t (z_{1}, z_{2}, \ldots, z_{m})$ and a ring $A$, 
let $A[\bm{z}] := A[z_{1}, z_{2}, \ldots, z_{m}]$ be the polynomial ring in these $m$ variables over $A$. 
We define two maps $L, R \colon G \to \Mat(n, \mathbb{C})$ as follows: 
\begin{align*}
L(g)_{ij} := 
\begin{cases}
1 & g = g_{i} g_{j}^{-1}, \\ 
0 & \text{otherwise},
\end{cases} \quad 
R(g)_{ij} := 
\begin{cases}
1 & g = g_{i}^{-1} g_{j}, \\ 
0 & \text{otherwise}. 
\end{cases}
\end{align*}
Note that $L$ and $R$ are matrix forms of the left and the right regular representations of $G$, respectively. 
Let $\bm{x} := {}^{t}(x_{g_{1}}, x_{g_{2}}, \ldots, x_{g_{n}})$. 
The group determinant $\Theta_{G}(x_{g})$ of $G$ is a homogeneous polynomial of degree $n$ defined by 
\[
\Theta_{G}(x_{g}) := \det{(x_{g h^{-1}})_{g, h \in G}} = \det{\left( \sum_{g \in G} x_{g} L(g) \right)} = \det{\left( \sum_{g \in G} x_{g} R(g) \right)} \in \mathbb{Z}[\bm{x}]. 
\]
Let $\partial \bm{x} := {}^t (\frac{\partial}{\partial x_{g_{1}}}, \frac{\partial}{\partial x_{g_{2}}}, \ldots, \frac{\partial}{\partial x_{g_{n}}})$ be the $n$-tuple of the partial differential operators. 
In Frobenius' study of the group determinant \cite[p.~50]{Frobenius1968gruppen}, 
the following partial differential operators were introduced: 
For each $g \in G$, 
\[
D_{g} := \sum_{h \in G} x_{g^{-1} h} \partial x_{h} = {}^t (L(g) \bm{x}) \partial \bm{x}, \quad 
\Delta_{g} := \sum_{h \in G} x_{h g^{-1}} \partial x_{h} = {}^t (R(g^{-1}) \bm{x}) \partial \bm{x}. 
\]
Let $\chi_{\varphi}$ denote the character of a representation $\varphi$ of $G$. 
Frobenius \cite[p.~50]{Frobenius1968gruppen} showed that for any irreducible representation $\varphi$ of $G$ over $\mathbb{C}$ and any $g \in G$,  
the following holds (see, e.g., \cite[p.~226]{MR803326}): 
\begin{align}
D_{g} \det{\left( \sum_{h \in G} \varphi(h) x_{h} \right)} 
= \Delta_{g} \det{\left( \sum_{h \in G} \varphi(h) x_{h} \right)} 
= \chi_{\varphi}(g) \det{\left( \sum_{h \in G} \varphi(h) x_{h} \right)}. 
\end{align}
Here, we may assume $\varphi$ to be an arbitrary representation of $G$ 
since every representation is equivalent to a direct sum of irreducible representations. 
Thus, in essence, Frobenius proved the following. 

\begin{thm}\label{thm:1}
Let $G$ be a finite group. 
For any representation $\varphi$ of $G$ over $\mathbb{C}$ and any $g \in G$, 
the following holds: 
\[
D_{g} \det{\left( \sum_{h \in G} \varphi(h) x_{h} \right)} 
= \Delta_{g} \det{\left( \sum_{h \in G} \varphi(h) x_{h} \right)} 
= \chi_{\varphi}(g) \det{\left( \sum_{h \in G} \varphi(h) x_{h} \right)}. 
\]
\end{thm}

From the definitions of $L$ and $R$, 
for any $g \in G \setminus \{ g_{1} \}$, 
we have $\chi_{L}(g) = \chi_{R}(g) = 0$. 
Hence, 
it follows from Theorem~$\ref{thm:1}$ that $D_{g} \Theta_{G}(x_{g}) = \Delta_{g} \Theta_{G}(x_{g}) = 0$ for any $g \in G \setminus \{ g_{1} \}$. 
This consideration leads to the study of the invariant rings 
\[
\mathbb{C}[\bm{x}]^{D_{g}} := \{ f(\bm{x}) \in \mathbb{C}[\bm{x}] \mid D_{g} f(\bm{x}) = 0 \}, 
\quad \mathbb{C}[\bm{x}]^{\Delta_{g}} := \{ f(\bm{x}) \in \mathbb{C}[\bm{x}] \mid \Delta_{g} f(\bm{x}) = 0 \}. 
\]
The first and second fundamental theorems of invariant theory concern the generators of an arbitrary invariant ring and the relations among these generators, where the first provides the generators and the second provides the relations. 
Recently, in the case that $n$ has at most two prime factors and $G$ is the cyclic group $\mathbb{Z} / n \mathbb{Z} = \left\{\overline{0}, \overline{1}, \ldots, \overline{n -  1} \right\}$, 
the first and second fundamental theorems of invariant theory for $\mathbb{C}[\bm{x}]^{D_{\overline{1}}}$ and $\mathbb{C}[\bm{x}]^{D_{\overline{n - 1}}}$ have been given \cite{Ochiai}.

Let $\widehat{G} := \{ \varphi_{1}, \varphi_{2}, \ldots, \varphi_{r} \}$ be a complete set of representatives of the equivalence classes of the irreducible representations of $G$, and let $\bm{y} := {}^{t}(y_{1}, y_{2}, \ldots, y_{r})$, where $y_{i} := \det{\left( \sum_{h \in G} \varphi_{i}(h) x_{h} \right)}$. 
In this paper, 
we provide a theorem concerning the invariant rings. 

\begin{thm}[]\label{thm:2}
Let $G$ be a finite group and let $S := \{ h_{1}, h_{2}, \ldots, h_{r} \}$ be a complete set of representatives of the conjugacy classes of $G$, where $h_{1} := g_{1}$. 
Then, the following holds: 
\[
\bigcap_{g \in S \setminus \{ h_{1} \} } \mathbb{C}[\bm{y}]^{D_{g}} 
= \bigcap_{g \in S \setminus \{ h_{1} \} } \mathbb{C}[\bm{y}]^{\Delta_{g}} 
= \mathbb{C}[\Theta_{G}(x_{g})]. 
\]
\end{thm}

We note that if $g, h \in G$ are conjugate in $G$, 
then, since every character is a class function, we have 
$D_{g} y_{i} = \Delta_{g} y_{i} = D_{h} y_{i} = \Delta_{h} y_{i}$ from Theorem~$\ref{thm:1}$, 
and therefore it holds that 
$\mathbb{C}[\bm{y}]^{D_{g}} = \mathbb{C}[\bm{y}]^{\Delta_{g}} = \mathbb{C}[\bm{y}]^{D_{h}} = \mathbb{C}[\bm{y}]^{\Delta_{h}}$. 

If $G$ is abelian, then it follows from the orthogonality relations that the matrix $P$ satisfying $\bm{y} = P \bm{x}$ is invertible. 
This means that $\mathbb{C}[\bm{y}] = \mathbb{C}[\bm{x}]$. 
Thus, as a corollary of Theorem~$\ref{thm:2}$, we obtain the following. 

\begin{cor}
Let $G$ be a finite abelian group. 
Then the following holds: 
\[
\bigcap_{g \in G \setminus \{ g_{1} \} } \mathbb{C}[\bm{x}]^{D_{g}} 
= \bigcap_{g \in G \setminus \{ g_{1} \} } \mathbb{C}[\bm{x}]^{\Delta_{g}} 
= \mathbb{C}[\Theta_{G}(x_{g})]. 
\]
\end{cor}

\begin{remark}
The following holds: 
Let $G$ be an abelian group of order $n$ and let 
\begin{align*}
\SL(\mathbb{C}G) &:= 
\left\{ 
(a_{g h^{-1}})_{g, h \in G} \mid a_{g_{1}}, a_{g_{2}}, \ldots, a_{g_{n}} \in \mathbb{C}
\right\} \cap \SL(n, \mathbb{C}), \\ 
\mathbb{C}[\bm{x}]^{\SL(\mathbb{C}G)} &:= \left\{ f(\bm{x}) \in \mathbb{C}[\bm{x}] \mid  f(A \bm{x}) = f(\bm{x}) \: \: \text{for any} \: A \in \SL(\mathbb{C}G) \right\}. 
\end{align*}
Then it holds that 
\[
\mathbb{C}[\bm{x}]^{\SL(\mathbb{C}G)} = \mathbb{C}[\Theta_{G}(x_{g})]. 
\]
This is a generalization of \cite[Theorem~4]{Ochiai}, which assumes that $G$ is a cyclic group. 
We can prove the equality in the same way as for the cyclic group. 
\end{remark}

In addition, in the case that $G$ is abelian, we give a necessary  and sufficient condition for the group determinant of a quotient group $G / H$ to be an element of $\mathbb{C}[\bm{x}]^{D_g}$. 

\begin{thm}\label{thm:5}
Let $G$ be a finite abelian group, 
let $H$ be a subgroup of $G$, and let 
\begin{align*}
\widehat{G}_{H} := \left\{ \chi \in \widehat{G} \mid \chi(h) = 1, \: h \in H \right\}. 
\end{align*}
Suppose that 
\begin{align*}
G = \displaystyle\bigsqcup_{t \in T} t H, \quad \widehat{G} = \displaystyle\bigsqcup_{\chi \in X} \chi \widehat{G}_{H}
\end{align*}
are residue class decompositions of $G$ and $\widehat{G}$, respectively. 
Then, for any $\chi \in X$, it holds that $D_{g} \Theta_{G / H}{\left( y_{t H}^{\chi} \right)} = 0$ if and only if $g \in G \setminus H$, 
where $y_{t H}^{\chi} := \sum_{h \in H} \chi(t h) x_{t h}$. 
\end{thm}

In relation to this theorem, 
we remark that the following is known as a generalized Dedekind's theorem \cite{MR4814714}: 
\begin{align*}
\Theta_{G}(x_{g}) = \prod_{\chi \in X} \Theta_{G / H}{\left( y_{t H}^{\chi} \right)} = \Theta_{H}(z_{h}), 
\end{align*}
where $z_{h} := \frac{1}{| H |}\sum_{\chi \in X} \chi(h^{-1}) \Theta_{G/H}{\left( y_{t H}^{\chi} \right)}$.

Also, we present the first and second fundamental theorems of invariant theory for $\mathbb{C}[\bm{y}]^{D_{g}}$ in the case that $G$ is the symmetric group $S_3$ of degree $3$. 

\begin{thm}[]\label{thm:6}
Let $G := S_3$, $H_1 := \langle (12) \rangle$, and $H_2 := \langle (123) \rangle$. Suppose that 
\[
G = t_1 H_1 \sqcup t_2 H_1 \sqcup t_3 H_1, \quad 
G = u_1 H_2 \sqcup u_2 H_2
\]
are residue class decompositions of $G$. 
Let $\widehat{H}_1 := \{\chi_1, \chi_2\}$ and $\widehat{H}_2 := \{\chi_3, \chi_4, \chi_5\}$ denote the character groups of $H_1$ and $H_2$, respectively, where $\chi_1$ and $\chi_3$ are the trivial characters. 
We define 
\begin{align*}
{\Theta'}_{G/H_1}^{\;\chi_k} := \det{\left( \sum_{h \in H_1} \chi_{k}(h) x_{t_{i} h t_{j}^{-1}} \right)_{1 \leq i, j \leq 3}}, \quad 
{\Theta'}_{G/H_2}^{\;\chi_l} := \det{\left( \sum_{h \in H_2} \chi_{l}(h) x_{u_{i} h u_{j}^{-1}} \right)_{1 \leq i, j \leq 2}} 
\end{align*}
for $k = 1, 2$ and $l = 3, 4, 5$. 
Then, the first and second fundamental theorems are given as follows: 
\[
\mathbb{C}[\bm{y}]^{D_{(123)}} = \mathbb{C}[{\Theta'}_{G/H_1}^{\; \chi_1}, {\Theta'}_{G/H_1}^{\;\chi_2}], \quad 
\mathbb{C}[\bm{y}]^{D_{(12)}} = \mathbb{C}[{\Theta'}_{G/H_2}^{\;\chi_3}, {\Theta'}_{G/H_2}^{\;\chi_4}]. 
\]
Here, any products of powers of ${\Theta'}_{G/H_1}^{\; \chi_1}$ and ${\Theta'}_{G/H_1}^{\; \chi_2}$ are linearly independent over $\mathbb{C}$, 
and any products of powers of ${\Theta'}_{G/H_2}^{\; \chi_3}$ and ${\Theta'}_{G/H_2}^{\; \chi_4}$ are linearly independent over $\mathbb{C}$. 
Also, the following relations hold: 
\[
\Theta_{G}(x_{g}) 
= \prod_{\chi \in \widehat{H}_{1}} {\Theta'}_{G/H_{1}}^{\;\chi} 
= \prod_{\chi \in \widehat{H}_{2}} {\Theta'}_{G/H_{2}}^{\;\chi}, \quad 
{\Theta'}_{G/H_2}^{\;\chi_4} = {\Theta'}_{G/H_2}^{\;\chi_5}. 
\]
\end{thm}

The structure of this paper is as follows. 
In Sections~\ref{sec:2}--\ref{sec:4}, we provide the proofs of Theorems~$\ref{thm:2}$, $\ref{thm:5}$, and $\ref{thm:6}$, respectively. 
In Section~\ref{sec:5}, we make a note on Theorem~\ref{thm:1}. Theorem~\ref{thm:1} was originally proved by Frobenius, 
but it should be noted that the definition of the character differs from the current one. 
With this difference in mind, we explain the Frobenius' proof of Theorem~\ref{thm:1}. 
In addition, we provide a more direct proof of Theorem~\ref{thm:1}.

\section{Proof of Theorem~$\mathbf{\ref{thm:2}}$}\label{sec:2}

To prove Theorem~$\ref{thm:2}$, we use the following facts. 
Let $d_{i} := \deg{\varphi_{i}}$ for $1 \leq i \leq r$. 
\begin{enumerate}
\item[(i)] Second orthogonality relations (see, e.g., \cite[p.~46, Theorem~4.4.12]{MR2867444}): Let $C$ and $C'$ be conjugacy classes of $G$, and let $g \in C$ and $h \in C'$. 
Then we have 
\begin{align*}
\sum_{i = 1}^{r} \chi_{\varphi_{i}}(g) \overline{\chi_{\varphi_{i}}(h)} = 
\begin{cases} 
\frac{|G|}{|C|} & C = C', \\ 
0 & C \neq C'. 
\end{cases}
\end{align*}
\item[(ii)] Frobenius determinant theorem \cite{Frobenius1968gruppen} (for a historical account, see \cite{MR1659232, MR1554141, MR803326}): We have 
\[
\Theta_{G} (x_{g}) = \prod_{i = 1}^{r} y_{i}^{d_{i}}. 
\]
\end{enumerate}

Below, we give a proof of Theorem~$\ref{thm:2}$. 
For a monomial $\bm{y^{\alpha}} := y_{1}^{\alpha_{1}} y_{2}^{\alpha_{2}} \cdots y_{r}^{\alpha_{r}}$ with $\bm{\alpha} := (\alpha_{1}, \alpha_{2}, \ldots, \alpha_{r}) \in (\mathbb{Z}_{\geq 0})^{r}$, it follows from Theorem~$\ref{thm:1}$ that $D_{g} \bm{y^{\alpha}} = \left(\sum_{i = 1}^{r} \chi_{\varphi_{i}}(g) \alpha_{i} \right)  \bm{y^{\alpha}}$ for any $g \in G$. 
To determine $\bm{\alpha} \in (\mathbb{Z}_{\geq 0})^{r}$ satisfying $D_{g} \bm{y^{\alpha}} = 0$ for any $g \in S \setminus \{h_{1} \}$, we consider the system of linear equations $M \bm{v} = \bm{0}$, 
where the coefficient matrix is given by  
\begin{align*}
M := 
\begin{pmatrix}
\chi_{\varphi_{1}}(h_{2}) & \chi_{\varphi_{2}}(h_{2}) & \cdots & \chi_{\varphi_{r}}(h_{2}) \\ 
\chi_{\varphi_{1}}(h_{3}) & \chi_{\varphi_{2}}(h_{3}) & \cdots & \chi_{\varphi_{r}}(h_{3}) \\ 
\vdots & \vdots & \ddots & \vdots \\ 
\chi_{\varphi_{1}}(h_{r}) & \chi_{\varphi_{2}}(h_{r}) & \cdots & \chi_{\varphi_{r}}(h_{r}) \\ 
\end{pmatrix}.
\end{align*}
From the second orthogonality relations, 
for any $c \in \mathbb{C}$, 
the vector 
\[
\bm{v} 
=
c 
\begin{pmatrix}
\chi_{\varphi_{1}}(h_{1}) \\ 
\chi_{\varphi_{2}}(h_{1}) \\ 
\vdots  \\ 
\chi_{\varphi_{r}}(h_{1}) 
\end{pmatrix} 
= 
c 
\begin{pmatrix}
d_{1} \\ 
d_{2} \\ 
\vdots  \\ 
d_{r} 
\end{pmatrix}
\]
is a solution of the system. 
Moreover, from the second orthogonality relations, we find that the row vectors of $M$ are linearly independent. 
Thus, the rank of $M$ is $r - 1$, so the dimension of the solution space of the system is $1$. 
Therefore, it holds that $D_{g} \bm{y^{\alpha}} = 0$ for any $g \in S \setminus \{h_{1} \}$ if and only if $\bm{y^{\alpha}} = y_{1}^{c d_{1}} y_{2}^{c d_{2}} \cdots y_{r}^{c d_{r}} = \Theta_{G} (x_{g})^{c}$ for some $c \in \mathbb{Z}_{\geq 0}$, where the last equality follows from the Frobenius determinant theorem.

\section{Proof of Theorem~$\mathbf{\ref{thm:5}}$}\label{sec:3}

We prove Theorem~$\ref{thm:5}$. 
Note that $\widehat{G}_{H} = \left\{ \chi' \circ \pi \: \vert \: \chi' \in \widehat{G/H} \right\}$ holds, where $\pi : G \rightarrow G/H$ denotes the canonical homomorphism. 
Let $y_{\chi} := \sum_{g \in G} \chi(g) x_{g}$ for $\chi \in \widehat{G}$. 
Then, for $\chi \in X$ and $\chi' \in \widehat{G}_{H}$, we have
\begin{align*}
y_{\chi \chi'} 
= \sum_{g \in G} (\chi \chi')(g) x_{g}  
= \sum_{t \in T} \sum_{h \in H} (\chi \chi')(t h) x_{t h}  
= \sum_{t \in T} \chi'(t) \sum_{h \in H} \chi(t h) x_{t h}  
= \sum_{t \in T} \chi'(t) y_{t H}^{\chi}, 
\end{align*}
and thus it holds that 
\begin{align*}
\prod_{\chi' \in \widehat{G}_{H}} y_{\chi \chi'} 
= \prod_{\chi' \in \widehat{G}_{H}} \sum_{t \in T} \chi'(t) y_{t H}^{\chi} 
= \prod_{\chi' \in \widehat{G/H}} \sum_{t H \in  G/H} \chi'(t H) y_{t H}^{\chi}  
= \Theta_{G/H}(y_{t H}^{\chi}), 
\end{align*}
where the last equality follows from the Frobenius determinant theorem. 
From this and Theorem~$\ref{thm:1}$, for any $\chi \in X$, we have 
\begin{align*}
D_{g} \Theta_{G/H}(y_{t H}^{\chi}) 
= D_{g} \prod_{\chi' \in  \widehat{G}_{H}} y_{\chi \chi'} 
= \sum_{\chi' \in \widehat{G}_{H}} (\chi \chi')(g) \prod_{\chi'' \in \widehat{G}_{H}} y_{\chi \chi''}  
= \Theta_{G/H}(y_{t H}^{\chi}) \chi(g) \sum_{\chi' \in \widehat{G}_{H}} \chi'(g). 
\end{align*}
Therefore, we find that 
\begin{align*}
D_{g} \Theta_{G/H}{\left( y_{t H}^{\chi} \right)} = 0 
\iff \!\!\sum_{\chi' \in \widehat{G}_{H}} \chi'(g) = 0 
\iff \!\!\sum_{\chi' \in \widehat{G/H}} \chi'(g H) \chi'(g_{1} H) = 0 
\iff g \in G \setminus H, 
\end{align*}
where the last equivalence follows from the second orthogonality relations.

\section{Proof of Theorem~$\mathbf{\ref{thm:6}}$}\label{sec:4}

We prove Theorem~$\ref{thm:6}$. 
From Theorem~$\ref{thm:1}$ and Table~$\ref{table1}$, for a monomial $\bm{y^{\alpha}} := y_{1}^{\alpha_{1}} y_{2}^{\alpha_{2}} y_{3}^{\alpha_{3}}$ with $\bm{\alpha} := (\alpha_{1}, \alpha_{2}, \alpha_{3}) \in (\mathbb{Z}_{\geq 0})^{3}$, 
we have $D_{(123)} \bm{y^{\alpha}} = (\alpha_{1} + \alpha_{2} - \alpha_{3}) \bm{y^{\alpha}}$ and $D_{(12)} \bm{y^{\alpha}} = (\alpha_{1} - \alpha_{2}) \bm{y^{\alpha}}$. 
Thus, it holds that $\mathbb{C}[\bm{y}]^{D_{(123)}} = \mathbb{C}[y_{1} y_{3}, y_{2} y_{3}]$ and $\mathbb{C}[\bm{y}]^{D_{(12)}} = \mathbb{C}[y_{1} y_{2}, y_{3}]$. 
Here, it follows from direct calculations that 
\begin{gather*}
y_{1} y_{3} = {\Theta'}_{G/H_{1}}^{\;\chi_{1}}, \quad 
y_{2} y_{3} = {\Theta'}_{G/H_{1}}^{\;\chi_{2}}, \quad 
y_{1} y_{2} = {\Theta'}_{G/H_{2}}^{\;\chi_3}, \quad 
y_{3} = {\Theta'}_{G/H_{2}}^{\;\chi_4} 
= {\Theta'}_{G/H_{2}}^{\;\chi_5}. 
\end{gather*}
Also, from the Frobenius determinant theorem, we have 
\begin{align*}
\Theta_{G}(x_{g}) 
= y_{1} y_{2} y_{3}^{2}
= \prod_{\chi \in \widehat{H}_{1}} {\Theta'}_{G/H_{1}}^{\;\chi} 
= \prod_{\chi \in \widehat{H}_{2}} {\Theta'}_{G/H_{2}}^{\;\chi}. 
\end{align*}
We show that any products of powers of $z_{1} := {\Theta'}_{G/H_1}^{\; \chi_1}$ and $z_{2} := {\Theta'}_{G/H_1}^{\; \chi_2}$ are linearly independent over $\mathbb{C}$. 
Since $z_{1}$ and $z_{2}$ are homogeneous polynomials in the $x_{g}$'s of the same degree, 
it is sufficient to consider the products of powers of $z_{1}$ and $z_{2}$ having a same degree. 
We prove by induction on $m$ that if $\sum_{i = 0}^{m} c_{i} z_{1}^{i} z_{2}^{m - i} = 0$ for some $c_{0}, c_{1}, \ldots, c_{m} \in \mathbb{C}$, then $c_{0} = c_{1} = \cdots = c_{m}$. 
The base case $m = 1$ is true since $z_{2}$ is not divisible by $z_{1}$. 
Suppose for induction that the case $m = k$ is true, and we now consider the case $m = k + 1$. 
Assume that $\sum_{i = 0}^{k + 1} c_{i} z_{1}^{i} z_{2}^{k + 1 - i} = 0$ for some $c_{0}, c_{1}, \ldots, c_{k + 1} \in \mathbb{C}$. 
Then, from $z_{2} \sum_{i = 0}^{k} c_{i} z_{1}^{i} z_{2}^{k - i} = - c_{k + 1} z_{1}^{k + 1}$, it follows that $c_{k + 1} = 0$ since $z_{2}$ is not divisible by $z_{1}$. 
This yields $\sum_{i = 0}^{k} c_{i} z_{1}^{i} z_{2}^{k - i} = 0$, 
and we have $c_{0} = c_{1} = \cdots = c_{k}$ by the inductive hypothesis. Thus, the case $m = k + 1$ is true. 
By mathematical induction, the above statement is true for any positive integer $m$. 
Therefore, any products of powers of $z_{1}$ and $z_{2}$ are linearly independent over $\mathbb{C}$. 
We can prove in the same way that any products of powers of ${\Theta'}_{G/H_2}^{\; \chi_3}$ and ${\Theta'}_{G/H_2}^{\; \chi_4}$ are linearly independent over $\mathbb{C}$. 


\begin{table}[H]
\centering
\[
\begin{array}{|c|cccccc|} \hline 
 & (1) & (12) & (13) & (23) & (123) & (132) \\
\hline
\varphi_{1} & 1 & 1 & 1 & 1 & 1 & 1 \\
\varphi_{2} & 1 & -1 & -1 & -1 & 1 & 1 \\
\varphi_{3} & 
\begin{pmatrix} 1 & 0 \\ 0 & 1 \end{pmatrix} & 
\begin{pmatrix} 0 & 1 \\ 1 & 0 \end{pmatrix} & 
\begin{pmatrix} 0 & \omega \\ \omega^2 & 0 \end{pmatrix} & 
\begin{pmatrix} 0 & \omega^2 \\ \omega & 0 \end{pmatrix} & 
\begin{pmatrix} \omega & 0 \\ 0 & \omega^2 \end{pmatrix} & \begin{pmatrix} \omega^2 & 0 \\ 0 & \omega \end{pmatrix} \\
\hline 
\end{array}
\]
\caption{Irreducible representations of $S_{3}$, where $\omega := \exp{\frac{2 \pi \sqrt{-1}}{3}}$.}
\label{table1}
\end{table}

\section{Note on Theorem~$\mathbf{\ref{thm:1}}$}\label{sec:5}

For any representation $\varphi$ of $G$, 
Frobenius \cite[pp.~44--45]{Frobenius1968gruppen} defined $\chi_{\varphi}$ as follows (see \cite[p.~166]{MR1554141} and \cite[p.~226]{MR803326} for a historical account):
\begin{align*}
\chi_{\varphi}(g) := 
\begin{cases}
d & g = g_{1}, \\ 
\text{the coefficient of} \: x_{g_{1}}^{d - 1} x_{g} \: \text{ in } \: \Phi := \det{\left( \sum_{h \in G} \varphi(h) x_{h} \right)} & g \neq g_{1},
\end{cases}
\end{align*}
where $d := \deg{\varphi}$. 
It is immediately follows from this definition that $\chi_{\varphi}(g)$ is equal to the value obtained by substituting $x_{h} = 1 \;(h = g_{1}); \, 0 \;(h \neq g_{1})$ into $\frac{\partial \Phi}{\partial x_{g}}$, so using this fact, Frobenius derived the equation~(1). 
Later, he \cite[p.~92]{Frobenius1897} defined $\chi_{\varphi}$ as the trace of $\varphi$. 
Since $\varphi(g_{1})$ is the identity matrix, we find that 
these two definitions of $\chi_{\varphi}$ are identical. 

Van der Waerden \cite[p.~226]{MR803326}, with respect to the equation~(1), states:  
\begin{quote}  
``As Hawkins notes, one can hardly overemphasize the importance of the relation.''  
\end{quote}  
He also mentions:  
\begin{quote}  
``Frobenius' derivation of the relation is very interesting. It has been explained by Th.~Hawkins in his paper \cite{MR526216}.''  
\end{quote}  
Hawkins \cite{MR0446837, MR526216} provides an explanation of Frobenius' derivation of the equation~(1), but it is not a reproduction of Frobenius' proof in \cite[p.~50]{Frobenius1968gruppen}.  
The Frobenius' proof is as follows (we remark that $\varphi$ does not need to be irreducible):  
Let $X := \sum_{h \in G} \varphi(h) x_{h}$, $Y := \sum_{h \in G} \varphi(h) y_{h}$, and $Z = \sum_{h \in G} \varphi(h) z_{h} := X Y$. 
Then, $z_{h} = \sum_{g \in G} x_{h g^{-1}} y_{g}$ holds. 
Thus, the chain rule for partial derivatives gives
\[
\frac{\partial \det{Z}}{\partial y_{g}} 
= \sum_{h \in G} \frac{\partial \det{Z}}{\partial z_{h}} \frac{\partial z_{h}}{\partial y_{g}} 
= \sum_{h \in G} \frac{\partial \det{Z}}{\partial z_{h}} x_{h g^{-1}}.
\]
On the other hand, 
it follows from $\det{Z} = \det{X} \det{Y}$ that 
\[
\frac{\partial \det{Z}}{\partial y_{g}} = \det{X} \frac{\partial \det{Y}}{\partial y_{g}}. 
\]
Therefore, we have 
\begin{align*}
\sum_{h \in G} \frac{\partial \det{Z}}{\partial z_{h}} x_{h g^{-1}} = \det{X} \frac{\partial \det{Y}}{\partial y_{g}}.
\end{align*}
Substituting $y_{h} = 1 \;(h = g_{1});\, 0 \;(h \neq g_{1})$, 
the left-hand side becomes $\sum_{h \in G} \frac{\partial \det{X}}{\partial x_{h}} x_{h g^{-1}}$ and the right-hand side becomes $\chi_{\varphi}(g) \det{X}$. 
From the above, we conclude that $\Delta_{g} \det{X} = \chi_{\varphi}(g) \det{X}$ holds. 
Also, by letting $Z := Y X$, 
we obtain $D_{g} \det{X} = \chi_{\varphi}(g) \det{X}$. 


Below, 
we provide a more direct proof of Theorem~$\ref{thm:1}$ by partitioning matrices. 
Let $X := \sum_{h \in G} \varphi(h) x_{h}$ and $\varphi(h) = (\varphi(h)_{i j})_{1 \leq i, j \leq d} = (\varphi^{(1)}(h) \;\; \varphi^{(2)}(h) \;\; \cdots \;\; \varphi^{(d)}(h))$. 
We define $\bm{v}_{j} := \sum_{h \in G} \varphi^{(j)}(h)x_{h}$ so that $X = \left( \bm{v}_{1} \;\; \bm{v}_{2} \;\; \cdots \;\; \bm{v}_{d} \right)$ holds. 
Then, it follows from $\frac{\partial}{\partial x_{h}} \det{X} 
= \sum_{j = 1}^{d} \det{\left( \bm{v}_{1} \;\; \cdots \;\; \bm{v}_{j - 1} \;\; \varphi^{(j)}(h) \;\; \bm{v}_{j + 1} \;\; \cdots \;\; \bm{v}_{d} \right)}$ that 
\begin{align*}
\Delta_{g} \det{X} 
= \sum_{h \in G} x_{h g^{-1}} \frac{\partial}{\partial x_{h}} \det{X} 
= \sum_{j = 1}^{d} \det{X_{j}},
\end{align*}
where $X_{j}$ denotes the matrix obtained by replacing the $j$-th row of $X$ with $\bm{v}'_{j} := \sum_{h \in G} \varphi^{(j)}(h) x_{h g^{-1}}$.  
Here, from 
\[
\sum_{h \in G} \varphi(h) x_{h g^{-1}} 
= \left( \sum_{h \in G} \varphi(h g^{-1}) x_{h g^{-1}} \right) \varphi(g) 
= \left( \sum_{h \in G} \varphi(h) x_{h} \right) \varphi(g),
\]
we have 
\begin{align*}
\sum_{h \in G} \varphi^{(j)}(h) x_{h g^{-1}}
= \left( \sum_{h \in G} \varphi(h) x_{h} \right) \varphi^{(j)}(g) 
= \left( \bm{v}_{1} \quad \bm{v}_{2} \quad \cdots \quad \bm{v}_{d} \right) 
\begin{pmatrix}
\varphi(g)_{1 j} \\ 
\varphi(g)_{2 j} \\ 
\vdots \\ 
\varphi(g)_{d j}
\end{pmatrix}. 
\end{align*}
That is, 
$\bm{v}'_{j} = \sum_{i = 1}^{d} \varphi(g)_{i j} \bm{v}_{i}$.
Thus, it holds that 
\begin{align*}
\det{X_{j}} 
&= \det{\left(\bm{v}_{1} \quad \cdots \quad \bm{v}_{j - 1} \quad \bm{v}'_{j} \quad \bm{v}_{j + 1} \quad \cdots \quad \bm{v}_{d} \right)} \\ 
&= \det{\left(\bm{v}_{1} \quad \cdots \quad \bm{v}_{j - 1} \quad \varphi(g)_{j j} \bm{v}_{j} \quad \bm{v}_{j + 1} \quad \cdots \quad \bm{v}_{d} \right)} \\ 
&= \varphi(g)_{j j} \det{X}.
\end{align*}
Therefore, we obtain 
\[
\Delta_{g} \det{X} = \sum_{j = 1}^{d} \varphi(g)_{j j} \det{X} = \chi_{\varphi}(g) \det{X}. 
\]
Also, by partitioning $X$ by rows, we can prove 
$D_{g} \det{X} = \chi_{\varphi}(g) \det{X}$.

\clearpage

\bibliography{reference}

\begin{thebibliography}{10}

\bibitem{MR1659232}
Keith Conrad.
\newblock The origin of representation theory.
\newblock {\em Enseign. Math. (2)}, 44(3-4):361--392, 1998.

\bibitem{Frobenius1968gruppen}
Ferdinand~Georg Frobenius.
\newblock \"{U}ber die {P}rimfactoren der {G}ruppendeterminante.
\newblock {\em Sitzungsberichte der K\"{o}niglich Preu{\ss}ischen Akademie der
  Wissenschaften zu Berlin}, pages 1343--1382, 1896.
\newblock Reprinted in {\it Gesammelte Abhandlungen, Band III}. Springer-Verlag
  Berlin Heidelberg, New York, 1968, pages 38--77.

\bibitem{Frobenius1897}
Ferdinand~Georg Frobenius.
\newblock \"{U}ber die {D}arstellung der endlichen {G}ruppen durch lineare
  {S}ubstitutionen.
\newblock {\em Sitzungsberichte der K\"{o}niglich Preu{\ss}ischen Akademie der
  Wissenschaften zu Berlin}, pages 944--1015, 1897.
\newblock Reprinted in {\it Gesammelte Abhandlungen, Band III}. Springer-Verlag
  Berlin Heidelberg, New York, 1967, pages 82--103.

\bibitem{MR1554141}
Thomas Hawkins.
\newblock The origins of the theory of group characters.
\newblock {\em Arch. History Exact Sci.}, 7(2):142--170, 1971.

\bibitem{MR0446837}
Thomas Hawkins.
\newblock New light on {F}robenius' creation of the theory of group characters.
\newblock {\em Arch. History Exact Sci.}, 12:217--243, 1974.

\bibitem{MR526216}
Thomas Hawkins.
\newblock The creation of the theory of group characters.
\newblock {\em Rice Univ. Stud.}, 64(2-3):57--71, 1978.
\newblock History of analysis (Proc. Conf., Rice Univ., Houston, Tex., 1977).

\bibitem{MR2867444}
Benjamin Steinberg.
\newblock {\em Representation theory of finite groups}.
\newblock Universitext. Springer, New York, 2012.
\newblock An introductory approach.

\bibitem{MR803326}
B.~L. van~der Waerden.
\newblock {\em A history of algebra}.
\newblock Springer-Verlag, Berlin, 1985.
\newblock From al-Khw\={a}rizm\={\i} to Emmy Noether.

\bibitem{Ochiai}
Naoya Yamaguchi, Hiroyuki Ochiai, and Yuka Yamaguchi.
\newblock First and second fundamental theorems for invariant rings generated
  by circulant determinants, 2025.
\newblock arXiv:2504.05688 [math.RT].

\bibitem{MR4814714}
Naoya Yamaguchi and Yuka Yamaguchi.
\newblock Generalized {D}edekind's theorem and its application to integer group
  determinants.
\newblock {\em J. Math. Soc. Japan}, 76(4):1123--1138, 2024.

\end{thebibliography}
\bibliographystyle{plain}

\medskip
\begin{flushleft}
Faculty of Education, 
University of Miyazaki, 
1-1 Gakuen Kibanadai-nishi, 
Miyazaki 889-2192, 
Japan \\ 
{\it Email address}, Yuka Yamaguchi: y-yamaguchi@cc.miyazaki-u.ac.jp \\ 
{\it Email address}, Naoya Yamaguchi: n-yamaguchi@cc.miyazaki-u.ac.jp 
\end{flushleft}

\end{document}